\documentclass[11pt]{article}
\usepackage{amsfonts, amsthm, amssymb, eucal, amsbsy}
\usepackage{amsmath}
\usepackage{chatterjee}

\begin{document}
\newtheorem{defn}{Definition}[section]
\newcommand{\gp}{g^\prime}
\newcommand{\gpp}{g^{\prime\prime}}
\newcommand{\gppp}{g^{\prime\prime\prime}}
\newcommand{\bs}{\boldsymbol{\sigma}}
\newcommand{\jj}{\mathcal{J}}

\title{A Simple Invariance Theorem}
\author{Sourav Chatterjee \\
U.C. Berkeley (Statistics) \\
sourav@stat.berkeley.edu}
\maketitle
\noindent {\it This is an old article (from May 2004), that will probably not be published, because a much improved paper with new results is in preparation. Still, I decided to put it in the archive because there are some things of interest here (in particular, the section on the S-K model) which will not appear in the new paper.}

\begin{abstract}
We present a simple extension of Lindeberg's 
argument for the Central Limit Theorem to get a general invariance
result. We apply the technique to prove results from random matrix
theory, spin glasses, and maxima of random fields.  
\end{abstract}

\section{Introduction and results}\label{results}
J.\ W.\ Lindeberg's elegant proof of the Central Limit Theorem
\cite{lindeberg20, lindeberg22}, despite being in the shadow
of Fourier analytic methods for a long time, is now
well known. It was revived by Trotter \cite{trotter59} and has since
been used successfully to derive CLTs in infinite dimensional spaces,
where the Fourier analytic methods are not so useful. For more
information on this topic, see the survey paper \cite{bgpr91} and the
monograph \cite{pr89}. (Another possible source is Bergstr\"om's books
\cite{bergstrom, bergstrom82}. It is also worth mentioning that LeCam
\cite{lecam} 
had a similar idea for Poisson approximation.) The ideas were carefully examined and generalized by Zolotarev \cite{zolotarev77} through the introduction of the so-called $\zeta$ metrics, which we shall not discuss here. 

However, it seems that the basic method of replacing non-Gaussian
random variables by Gaussians one by one and using Taylor expansion to
get approximation bounds has been applied only for proving central
limit theorems for sums of independent random elements, and its
potential for proving more general invariance results has been
overlooked in the literature. (After the preparation of the initial draft of this article, it came to our notice that indeed, there is an old article of Rotar \cite{rotar79} which examines the Lindeberg method polynomial maps in a limiting case. Also, earlier this year, Mossel, O'Donnell, and Oleszkiewicz \cite{mossel05} made some striking applications to problems from computer science and discrete mathematics using the Lindeberg method on polynomials.)

We shall derive a very simple extension
of Lindeberg's argument to obtain a result for general smooth functions. Basically, we shall
show that if $f:\rr^n \ra \rr$ is a function such that reasonable
fluctuations in any single coordinate (keeping others fixed) do not
affect the value of the function in a ``big'' way, then the
distribution of $f(X_1,\ldots,X_n)$, where $X_i$'s are independent
random variables, depends mainly on the first two moments of the
$X_i$'s. 

To make things precise, we first need a suitable measure of the
largest possible influence of any single coordinate on the outcome. 
\begin{defn}\label{lamdef}
For any open interval $I$ containing $0$, any positive integer $n$,
any function $f:I^n \ra \cc$ which is thrice differentiable in each
coordinate, and $1\le r\le 3$, let
\[
\lambda_r(f) := \sup\{|\partial_i^p f(\bx)|^{r/p}: 1\le i \le n,\ 1\le p
\le r, \ \bx \in I^n\}
\]
where $\partial_i^p$ denotes $p$-fold differentiation with respect to
the $i^{\mathrm{th}}$ coordinate. For a collection $\mf$ of such
functions, define $\lambda_r(\mf) := \sup_{f\in \mf} \lambda_r(f)$. 
\end{defn}
Note that the interval $I$ can be bounded or unbounded.
The numbers $\lambda_r(f)$ jointly constitute a measure of the maximum
possible influence of the fluctuation in a single coordinate on the
value of $f$ at any point in the set $I^n$. We shall show that $f$
will have the aforementioned invariance property when $\lambda_2(f)$
and $\lambda_3(f)$ are sufficiently small. 

In this paper, we shall generally denote vectors by $\bx, \by$
etc. The $i^{\mathrm{th}}$ component of $\bx$ will be denoted by
$x_i$, of $\by$ by $y_i$ and so on. 

In what follows, $\bbx = (X_1,\ldots,X_n)$ and $\bby
=(Y_1,\ldots,Y_n)$ are two vectors of independent random variables
with finite second moments, taking values in some open interval $I$ and 
satisfying, for each $i$, $\ee X_i = \ee
Y_i$ and $\ee X_i^2 = \ee Y_i^2$. We shall also assume that $\bbx$ and
$\bby$ are defined on the same probability space and are
independent. Finally, let $\gamma = \max\{\ee|X_i|^3, \ee|Y_i|^3, 1\le 
i\le n\}$. Note that $\gamma$ may be $\infty$.\\
\\
Here is our main result:
\begin{thm}\label{boundthm}
Let $f:I^n \ra\rr$ be thrice
differentiable in each argument. If we set $U = f(\bbx)$ and
$V = f(\bby)$, then for any thrice differentiable $g:\rr \ra \rr$ and
any $K > 0$,   
\begin{small}
\begin{eqnarray*}
|\ee g(U) - \ee g(V)| &\le& C_1(g) \lambda_2(f) \sum_{i=1}^n [\ee
(X_i^2; |X_i| > K) + \ee (Y_i^2; |Y_i| > K)] \\
& & + \ C_2(g) \lambda_3(f) \sum_{i=1}^n [\ee (|X_i|^3; |X_i| \le K) +
 \ee (|Y_i|^3; |Y_i| \le K)] 
\end{eqnarray*}
\end{small}
where $C_1(g) = \|\gp\|_\infty + \|\gpp\|_\infty$ and
$C_2(g) = \frac{1}{6}\|\gp\|_\infty + \frac{1}{2}\|\gpp\|_\infty +
\frac{1}{6}\|\gppp\|_\infty$.
\end{thm}
The last term in the above bound is usually dealt with as follows:
having chosen a suitable $K$, we use $\ee (|X_{i}|^3; |X_{i}|\le K)
\le K\ee (X_{i}^2)$. When $\gamma <\infty$, we can do better: 
\begin{cor}\label{cor1}
In the setting of the above Theorem, if we further have $\gamma
<\infty$, then $|\ee g(U) - \ee g(V)| \le 2C_2(g) \gamma
n\lambda_3(f)$.  
\end{cor}
For a quick example to see how Theorem \ref{boundthm} can be applied,
consider the function $f(\bx) = n^{-1/2}\sum_{i=1}^n x_i$. It is very easy to
compute $\lambda_2(f) = n^{-1}$ and $\lambda_3(f) = n^{-3/2}$. Now
suppose $X_i$'s are i.i.d.\ and 
$Y_i$'s are also i.i.d. Further, assume $\ee X_i = \ee Y_i = 0$
and $\ee X_i^2 = \ee Y_i^2 = 1$ for all $i$. Then taking $K = \epsilon 
\sqrt{n}$ and using Theorem \ref{boundthm} we can easily get
\begin{eqnarray*}
|\ee g(\frac{1}{\sqrt{n}}\sum_{i=1}^n X_i) - \ee
g(\frac{1}{\sqrt{n}}\sum_{i=1}^n Y_i)| 
&\le& C_1(g) [\ee(X_1^2; |X_1| > \epsilon \sqrt{n}) \\
& & + \ \ee(Y_1^2; |Y_1| > \epsilon \sqrt{n})] + 2C_2(g) \epsilon.
\end{eqnarray*}
Taking $n \ra \infty$, this proves the classical CLT since $\epsilon$
is arbitrary. Furthermore, if we assume that $\ee|X_1|^3 <\infty$ and
$\ee |Y_1|^3 < \infty$, then we also get an explicit error bound: 
\[
|\ee g(\frac{1}{\sqrt{n}}\sum_{i=1}^n X_i) - \ee
g(\frac{1}{\sqrt{n}}\sum_{i=1}^n Y_i)| \le \frac{C_2(g)[\ee|X_1|^3 +
  \ee |Y_1|^3]}{\sqrt{n}}.
\] 
For a more complicated example, consider the Stieltjes transform of a Wigner 
matrix. For a given $z\in \cc \backslash \rr$, define a function $f$
as 
\[
f((x_{ij})_{1\le i\le j\le N}) = \frac{1}{N} \tr ((A((x_{ij})) - zI)^{-1})
\]
where $A((x_{ij}))$ is the $N$ by $N$ matrix whose $(i,j)^{\mathrm{th}}$
element is $N^{-1/2}x_{ij}$ if $i\le j$ and
$N^{-1/2}x_{ji}$ otherwise, and $I$ is the $N$ by $N$ 
identity matrix, and ``tr'' stands for the trace of a matrix. In section
\ref{matrix} we shall use Theorem \ref{boundthm} to obtain invariance
results about this function, which will in turn yield the
weakest known condition for convergence of spectral measures to
Wigner's semicircle law. 

Another nontrivial example that we shall consider (in section
\ref{spinglass}) is the {\it free
  energy} of the Sherrington-Kirkpatrick model of spin glass theory.
Here the function $f$ is given by
\[
f((x_{ij})_{1\le i< j\le N}) = \frac{1}{N} \log \biggl[\sum_{\bs}
  \exp\bigl\{\frac{\beta}{\sqrt{N}} \sum_{i<j} x_{ij}\sigma_i \sigma_j + \beta
    h \sum_i \sigma_i\bigr\}\biggr]
\]
where the sum is taken over all $\bs = (\sigma_1,\ldots, \sigma_N) \in 
\{-1,1\}^N$, and $\beta, h$ are parameters. To deal with functions of
this form, which commonly occur as free energy functions of various
physical models, we have the following general Theorem:
\begin{thm}\label{likelihood}
Suppose $\mf$ is a finite collection of coordinatewise thrice 
differentiable functions from $I^n$ into $\rr$, and $\alpha \ge 1$. If
$F:I^n \ra \rr$ is defined as $F(\bx) := \alpha^{-1}\log[\sum_{f\in
  \mf} e^{\alpha f(\bx)}]$, then $\lambda_2(F)\le 3\alpha
\lambda_2(\mf)$ and $\lambda_3(F) \le 13 \alpha^2\lambda_3(\mf)$. 
\end{thm} 
In section \ref{spinglass}, we shall derive a condition under which
the asymptotic behaviour of the free energy in the
Sherrington-Kirkpatrick model is not dependent on the exact
distributions of the entries. Our condition
is weaker than the weakest known condition. In particular, it includes 
the ``i.i.d.\ mean zero unit variance'' case.

Besides the possible applications to free energy functions as mentioned
before, Theorem \ref{likelihood} can have other important uses, as
well. For example, the following result is an easy application of
Theorems \ref{boundthm} and \ref{likelihood}:  
\begin{thm}\label{maxthm}
Let $\mf$ be as in Theorem \ref{likelihood}. Let $U = \max_{f\in \mf}
f(\bbx)$ and $V = \max_{f\in \mf} f(\bby)$. Then for any thrice
differentiable $g:\rr \ra \rr$, any $K >0$, and any $\alpha \ge 1$, we 
have
\begin{eqnarray*}
|\ee g(U) - \ee g(V)|
&\le& 2\|\gp\|_\infty\alpha^{-1} \log|\mf| + 3\alpha C_1(g)
\lambda_2(\mf) T_1(K) \\ 
& & \ + \ 13\alpha^2 C_2(g)
\lambda_3(\mf) T_2(K)
\end{eqnarray*}
where $T_1(K) = \sum_{i=1}^n [\ee (X_i^2; |X_i| > K) + \ee (Y_i^2;
|Y_i| > K)]$ and $T_2(K) = \sum_{i=1}^n [\ee (|X_i|^3; |X_i| \le K) +
\ee (|Y_i|^3; |Y_i| \le K)]$.
\end{thm}
Again, we shall usually deal with $T_2(K)$ using $\ee (|X|^3; |X|\le
K) \le K \ee X^2$. If $\gamma < \infty$, we have a more explicit bound:
\begin{cor}\label{cor2}
In the setting of the above Theorem, if we further have $\gamma
<\infty$, then 
\[
|\ee g(U) - \ee g(V)| \le
K(g) [(\gamma n\lambda_3(\mf))^{1/3} (\log|\mf|)^{2/3} + \gamma
n\lambda_3(\mf)] 
\]
where $K(g) = \frac{19}{3}\|\gp\|_\infty + 13 \|\gpp\|_\infty +
\frac{13}{3}\|\gppp\|_\infty$.
\end{cor}
In section \ref{ground}, we shall demonstrate an application of
Theorem \ref{maxthm} involving the energy of the ground state in the
Sherrington-Kirkpatrick model of spin glasses. Essentially, we shall
show that under the same conditions on the $x_{ij}$'s as in section
\ref{spinglass}, the asymptotic behaviour of
\[
N^{-3/2}\max_{\bs} \sum_{1\le i<j\le N} x_{ij} \sigma_i\sigma_j,
\]
where the maximum is taken over all $\bs\in \{-1,1\}^N$, is not
dependent on the exact distributions of the $x_{ij}$'s. 

For an immediate application, consider the (very old) question
raised by Erd\H{o}s and Kac \cite{erdoskac46}: what is the limiting
distribution 
of $\max_{1\le j\le n} \frac{1}{\sqrt{n}} \sum_{i=1}^j X_i$ 
where $X_i$'s are i.i.d.\ with mean zero and unit variance? It is now
well known that the limiting distribution is the same as that of
$|Z|$, where $Z\sim N(0,1)$. Erd\H{o}s and Kac proved the result for the 
case of the simple random walk; the general result could be proved
only after Donsker established the weak invariance principle. Using
Corollary \ref{cor2}, we can easily establish concrete error bounds under
finite third moments assumption for this problem.

To work things out, let $\mf = \{f_i: 1\le i\le n\}$, where $f_i(\bx)
:= n^{-1/2} \sum_{j=1}^i x_j$. Clearly, $\lambda_3(\mf) = \max_{1\le i
  \le n} \lambda_3(f_i) = n^{-3/2}$. Corollary \ref{cor2} now  gives
the bound 
\[
|\ee g(U) - \ee g(V)| \le K(g) [\gamma^{1/3} n^{-1/6}
(\log n)^{2/3} + \gamma n^{-1/2}]
\]
where $U= \max_{1\le i\le n} \frac{1}{\sqrt{n}}\sum_{j=1}^i X_j$ and
$V= \max_{1\le i\le n} \frac{1}{\sqrt{n}}\sum_{j=1}^i Y_j$.

The three Theorems presented in this section are very general in
applicability, and present a unifying approach to solving examples
of the kind mentioned above, rather than applying different techniques 
for different problems. However, the method has its deficiencies, the
greatest being that functions have to be smooth. This is a rather
severe restriction, 
and eliminates a lot of interesting examples. For example, the method
will not allow us to deal with non-smooth functionals like stopping
times (in the case of random walks) and empirical distribution
functions (for random matrices). Smoothing approximations may
sometimes give crude bounds. Furthermore, the restriction about the
boundedness of derivatives hampers the applicability to many
interesting functions like spectral radii of random matrices. Again,
truncation techniques might work.  

The next three sections will be devoted to working out in detail the
examples mentioned before. Proofs of the Theorems and Corollaries will
be presented in the last section.

\section{Convergence of spectral distributions}\label{matrix}
In this section, we shall illustrate the application of our method to
proving invariance results about random matrices. Specifically, we
shall derive the weakest known condition under which the spectral
measures of a sequence of Wigner matrices converge to the semicircle
law. We begin with a very short introduction to some material from the 
spectral theory of large dimensional random matrices.

\subsection{Spectral measures}
The {\it Empirical Spectral Distribution} (ESD) of a square matrix 
is the probability distribution on the complex plane which puts
equal mass on each eigenvalue of the matrix (repeated by
multiplicities). The limit of a
sequence of ESDs is called the {\it Limiting Spectral Distribution}
(LSD) of the corresponding sequence of matrices. The existence and
identification of LSDs for various kinds of random matrices is one of
the main goals of random matrix theory.

For an excellent review of mathematical results known about limiting
spectral behaviour and further references, see Bai \cite{bai99}. For
relevance in physics, see the book by Mehta \cite{mehta91}.

\subsection{Stieltjes transforms}
A standard tool for identifying the LSD of a sequence of random
matrices is the Stieltjes transform. To cut a long story short, we can 
say that the ESDs of a sequence $\{A_N\}_{N=1}^\infty$ of random real
symmetric 
matrices converge in probability (w.r.t.\ the Prokhorov metric, for
example) to a probability distribution $G$ if and only if  
\[
\forall z \in \cc\backslash \rr, \ \frac{1}{N} \tr((A_N - zI_N)^{-1})
\stackrel{P}{\longrightarrow} \int_{-\infty}^\infty \frac{1}{x-z} dG(x)
\]
where $I_N$ is the identity matrix of order $N$. The expression on the 
right is the Stieltjes transform of $G$ evaluated at $z$. Similarly,
the expression on the left is the Stieltjes transform of the ESD of
$A_N$, evaluated at $z$. Stieltjes transforms
will be particularly useful for applying our technique, since they are 
infinitely differentiable as functions of the matrix entries.

\subsection{Wigner matrices}
A random Wigner matrix of order $N$ is an $N$ by $N$ real symmetric
matrix with independent entries on and above the diagonal.

More specifically, consider the map $A$ which ``constructs''
Wigner matrices of order $N$. Let $n= N(N+1)/2$ and write elements of
$\rr^n$ as $\bx = (x_{ij})_{1\le i\le j\le N}$. For any $\bx \in
\rr^n$, let $A(\bx)$ be the matrix whose 
$(i,j)^{\mathrm{th}}$ entry is $N^{-1/2}x_{ij}$ if $i\le j$ and
$N^{-1/2}x_{ji}$ if $i>j$. If $\bbx$ is a vector of $n$ independent
standard Gaussian random variables, then $A(\bbx)$ is a standard Gaussian
Wigner matrix. Wigner \cite{wigner55} showed that the LSD for a
sequence of standard Gaussian Wigner matrices is the semicircle law,
which has density $(2\pi)^{-1}\sqrt{4-x^2}$ in $[-2,2]$. 

It was later shown that the distribution of the 
entries do not play a significant role: convergence to the
semicircle law would hold under more general 
conditions (Cf. Arnold \cite{arnold67}, Grenander
\cite{grenander63} and Bai \cite{bai99}). The weakest known
condition under which the convergence to semicircle law holds was
given by Pastur \cite{pastur72}. It is claimed that the condition was 
shown to be necessary by Girko \cite{girko88}. For a
detailed exposition, see Bai \cite{bai99} or Khorunzhy, Khoruzhenko
and Pastur \cite{kkp96}. 

The method of this paper will give an easy way to show the sufficiency
of Pastur's condition. Incidentally, somewhat similar ideas involving
derivatives of empirical characteristic functions (instead of
Stieltjes transforms) to get concenration bounds for ESDs have been
explored in Chatterjee and Bose \cite{chatterjeebose}. 

\subsection{Derivation of Pastur's condition}
To get started, fix $z = u + iv \in \cc$, with $v \ne 0$. Define
$f:\rr^n \ra \rr$ as  
\[
f(\bx) := \frac{1}{N} \tr((A(\bx) - zI)^{-1}).
\]
Also, define $G : \rr^n \ra \cc^{N\times N}$ as $G(\bx) := (A(\bx) -
zI)^{-1}$. Now note that from matrix
theory we know that inverting a matrix involves computing the
classical adjoint and dividing by the determinant, which implies that
the elements of the inverse are all rational functions of the elements 
of the original matrix. Also note that since all eigenvalues of
$A(\bx)$ are real, therefore $\det(A(\bx)-zI) \ne 0$. Thus, $G$ is
infinitely differentiable along each coordinate. Also note that
$(A(\bx)-zI)G(\bx) = I$ for each $\bx$. Thus for $1\le i \le j \le N$,
$\frac{\partial}{\partial x_{ij}} [(A-zI)G] \equiv 0$, 
which gives 
\[
\frac{\partial G}{\partial x_{ij}} = -G\frac{\partial A}{\partial x_{ij}} G.
\]
Also, note that higher order derivatives of $A$ vanish
identically. Combining everything we easily get
\begin{eqnarray}
\frac{\partial f}{\partial x_{ij}} &=& -\frac{1}{N}\tr(\frac{\partial
  A}{\partial x_{ij}} G^2), \label{der1}\\ 
\frac{\partial^2 f}{\partial x_{ij}^2} &=& \frac{2}{N} \tr(\frac{\partial
  A}{\partial x_{ij}} G \frac{\partial  A}{\partial x_{ij}} G^2),
\label{der2}\\ 
\frac{\partial^3 f}{\partial x_{ij}^3} &=& -\frac{6}{N} \tr (\frac{\partial
  A}{\partial x_{ij}} G\frac{\partial  A}{\partial x_{ij}} G\frac{\partial
  A}{\partial x_{ij}}G^2 ). \label{der3}
\end{eqnarray}
Now we need to find good bounds for the above quantities. For that, we 
need some preparation. 

For an $N\times N$ complex matrix $B = ((b_{ij}))$, the
\emph{Hilbert-Schmidt norm} (or Schur norm, or Euclidean norm) of $B$
is defined as $\|B\| := (\sum_{i,j} |b_{ij}|^2)^{1/2}$. Besides the
usual properties of a matrix norm, it also satisfies the following:
\begin{enumerate}
\item $|\tr(BC)| \le \|B\|\|C\|$.
\item If $U$ is a unitary matrix, then for any $C$ of the same order,
  $\|CU\| = \|UC\| =  \|C\|$.
\item For a normal matrix $B$ (i.e.\  $B^*B = BB^*$, $B^*$ being the
  conjugate transpose of $B$) with eigenvalues 
  $\lambda_1,\ldots \lambda_N$, and any $C$, $\max\{\|BC\|,\|CB\|\}
  \le \max_{1\le i \le N} |\lambda_i|\cdot \|C\|$.
\end{enumerate}
The first property follows from the Cauchy-Schwarz inequality. The
second is true because $\|U\by\|_2 = \|\by\|_2$ for any
unitary matrix $U$ and any vector $\by \in \rr^N$, where $\|\cdot\|_2$
denotes the Euclidean norm on $\rr^N$.
For the last one, note that any normal matrix $B$ can be written as $B
= U\Delta U^*$ where $U$ is unitary and $\Delta$ is diagonal, with the
diagonal elements being the eigenvalues of $B$, and then apply the
second property.

The above facts are standard, and may be looked up in any standard
text on matrix analysis. See Wilkinson \cite{wilkinson67} pp.\ 55--58,
for example. 

Now, it is easy to see that $G$ and the derivatives of $A$ are all
normal matrices. Moreover, the eigenvalues of $G$ are bounded by
$|v|^{-1}$ (where $v = \mathrm{Im} \ z$) and the 
eigenvalues of $\partial A/ \partial x_{ij}$ are bounded by
$N^{-1/2}$. (Note that $\partial A / \partial x_{ij}$ is the matrix
which has $N^{-1/2}$ at the $(i,j)^{\mathrm{th}}$ and
$(j,i)^{\mathrm{th}}$ positions, and $0$ elsewhere.)

Thus, from the spectral representation of $G^2$ it follows that the
elements of $G^2$ are bounded by $|v|^{-2}$. This fact, and the identity 
(\ref{der1}) imply that
\begin{eqnarray}\label{bd1}
\left\|\frac{\partial f}{\partial x_{ij}}\right\|_\infty \le 2|v|^{-2}
N^{-3/2}.
\end{eqnarray}
Next, using the expression (\ref{der2}) and the three properties of
the Hilbert-Schmidt norm discussed above, we get 
\begin{eqnarray}\label{bd2}
\left\|\frac{\partial^2 f}{\partial x_{ij}^2}\right\|_\infty &\le&
\frac{2}{N} \left\|\frac{\partial A}{\partial x_{ij}}\right\| \left\|G 
  \frac{\partial A}{\partial x_{ij}}G^2\right\| \ \le \
4|v|^{-3}N^{-2}.
\end{eqnarray}
Similarly, (\ref{der3}) gives
\begin{eqnarray}\label{bd3} 
\left\|\frac{\partial^3 f}{\partial x_{ij}^3}\right\|_\infty &\le& 12
|v|^{-4} N^{-5/2}.
\end{eqnarray}
From (\ref{bd1}), (\ref{bd2}) and (\ref{bd3}) it follows that
\begin{eqnarray*}
\lambda_2(f) &\le& 4 \max\{|v|^{-4}, |v|^{-3}\} N^{-2},  \\
\lambda_3(f) &\le& 12 \max\{|v|^{-6}, |v|^{-4}\} N^{-5/2}.
\end{eqnarray*}
Let $\bbx = (X_{ij})_{1\le i\le j\le N}$ and $\bby = (Y_{ij})_{1\le i
  \le j \le N}$ be collections of independent random variables with
zero mean and unit variance. Let $U = \mathrm{Re}\, f(\bbx)$ and $V=
\mathrm{Re} \, f(\bby)$, and let $g:\rr \ra \rr$ be any thrice 
differentiable function. Note that $\mathrm{Re}\, f$ is a smooth
function and $\lambda_r(\mathrm{Re}\, f) \le \lambda_r(f)$ for each
$r$. With $K = \epsilon \sqrt{N}$, Theorem
\ref{boundthm} immediately tells us that $|\ee g(U) - \ee g(V)|$ can
be bounded by a multiple (depending only on $g$ and $v$) of
\begin{eqnarray*}
N^{-2} \sum_{1\le
  i \le j\le N} [\ee (X_{ij}^2; |X_{ij}| > \epsilon \sqrt{N}) +
\ee(Y_{ij}^2; |Y_{ij}| > \epsilon \sqrt{N})] + \epsilon.
\end{eqnarray*}
The same bound also works for functions of the imaginary parts. Using
this result and Wigner's Theorem for Gaussian matrices, we see 
that convergence to the semicircle law holds whenever $X_{ij}$'s are
independent with zero mean and unit variance, and satisfy
\begin{eqnarray}\label{mincond} 
& & \forall \epsilon > 0, \ \lim_{N\ra \infty}  N^{-2}\sum_{1\le  i \le
  j\le N} \ee (X_{ij}^2; |X_{ij}| > \epsilon \sqrt{N}) = 0.  
\end{eqnarray}
This is exactly Pastur's condition, as mentioned before. The condition 
is satisfied, for example, if $X_{ij}$'s are i.i.d.\ with zero mean
and unit variance. Also note that though this looks like Lindeberg's
condition for the central limit theorem, it is not exactly that.

\section{Universality of a spin glass model}\label{spinglass}
In this section, we obtain a condition for invariance (or, as
physicists say, universality) of the limiting free energy of the
Sherrington-Kirkpatrick model of spin glasses. We begin with a short
introduction. 

\subsection{Spin glasses}
Let $\Sigma_N = \{-1,1\}^N$. This is the space of all possible {\it
  spins} of $N$ particles in statistical mechanics. The spins are
random, but not independent --- the spin of one particle exerts
influence on the spin of another. The joint law of the $N$ spins is a
matter of great interest and intrigue. Various models have been
suggested over the years for various situations. Some of these models, 
like the famous {\it Ising} model, are deterministic in the sense that 
none of the model parameters are random, while some others, 
like the {\it Sherrington-Kirkpatrick} model which we shall discuss
here, involve random variables as model parameters.

All models assign a probability proportional to $\exp(-\beta
H_N(\bs))$ to the configuration $\bs$, where $H_N$ is the {\it
  Hamiltonian}, and $\beta = 1/T$, $T$ being the {\it
  temperature}. The {\it partition function} is $Z_N = \sum_{\bs}
\exp(-\beta H_N(\bs))$, and the {\it free energy} is the log of the
partition function divided by $N$. The asymptotic behaviour of the free
energy is of great consequence and interest to physicists, and
nowadays, to people in neural networks also.

For a detailed discussion of mathematical results about spin glass
models and further references, see Talagrand \cite{talagrand03}, for
instance. 

\subsection{The Sherrington-Kirkpatrick model}
The Sherrington-Kirkpatrick (S-K) model, introduced in \cite{sk72},
can be briefly described as follows: For each $N\ge 1$ 
let $\{\jj^N_{ij}, \ 1\le i, j \le N\}$ be a  
collection of i.i.d.\ $N(0,1)$ random variables. The S-K model assigns  
a random probability distribution (the Gibbs measure) on $\Sigma_N$ as
follows: For any configuration $\bs \in \Sigma_N$, the probability of
the system being in the state $\bs = (\sigma_1,\ldots,\sigma_N)$ is
given by 
\[
p_{N, \jj}(\bs) = Z_{N,\jj}^{-1} \exp(-\beta H_{N,\jj}(\bs))
\]
where $H_{N,\jj}(\bs) = - \frac{1}{\sqrt{N}} \sum_{i<j} \jj^N_{ij} \sigma_i
\sigma_j - h\sum_{i\le N} \sigma_i$, $\beta$ and $h$ are fixed
parameters, and $Z_{N,\jj}$ is the normalising
constant. Ideally, the subscripts should include $\beta$ and $h$, but
we are considering them to be fixed. It has been shown by Guerra and Toninelli
\cite{guerratoninelli02} that the limit 
\[
\lim_{N\ra \infty} \frac{1}{N} \ee(\log Z_{N,\jj})
\]
exists for all $\beta$ and $h$. See Talagrand \cite{talagrand03}
Theorem 2.10.1, p.\ 140 for a proof. A formula for the limit 
was conjectured by Parisi and proved by Talagrand
\cite{talagrand03b}. Talagrand (\cite{talagrand03} Corollary 2.2.5,
p.\ 32) also proves (in particular) that 
\[
\frac{1}{N} (\log Z_{N,\jj} - \ee \log Z_{N,\jj}) \stackrel{P}{\ra} 0
\]
for any $\beta$ and $h$. Both the above facts were proved under the
condition that $\jj^N_{ij}$ are  
i.i.d.\ $N(0,1)$. In fact, the rigorous proofs involve the use of
intricate properties of Gaussian random variables. Recently, in a
paper which was archived at a time when this article was being
written, Carmona and Hu \cite{carmonahu04} have proved that the limit
will exist and be the same when $\jj_{ij}$ are i.i.d.\ with zero mean, unit
variance and finite third moment. Their technique may be extended to
the case of independent variables with uniformly bounded third
absolute moments.  

We shall derive a sufficient condition for invariance of the limiting
free energy, which is weaker than the condition given by Carmona and
Hu, and includes the case where $\jj_{ij}$'s are i.i.d.\ with 
zero mean and unit variance, with no assumption about the third moment.

\subsection{Our condition}
Let $\mf = \{f_{\bs}: \bs \in \{-1,1\}^N\}$, where 
\[
f_{\bs}((x_{ij})) = \beta N^{-3/2} \sum_{i<j} x_{ij} \sigma_i\sigma_j
+ \beta h N^{-1}\sum_i \sigma_i.
\]
Then clearly, $\lambda_2(\mf) =
\beta^2 N^{-3}$, $\lambda_3(\mf) = \beta^3 N^{-9/2}$ and 
$|\mf| = 2^N$. Now, if we define $F(\bx) = N^{-1} \log [\sum_{\bs} 
e^{Nf_{\bs}(\bx)}]$, then by Theorem \ref{likelihood}, $\lambda_2(F)
\le 3\beta^2N^{-2}$ and $\lambda_3(F) \le 13\beta^3N^{-5/2}$. 

Suppose $\jj$ and $\jj^\prime$ are collections of independent random 
variables with zero mean and unit variance. If we let $U_N = F(\jj)$
and $V_N = F(\jj^\prime)$, then by Theorem 
\ref{boundthm}, for any thrice differentiable $g:\rr \ra \rr$ and any
fixed $\epsilon>0$, $|\ee g(U_N) - \ee g(V_N)|$ is bounded by
a constant multiple (depending only on $g$ and $\beta$) of 
\begin{eqnarray*}
N^{-2} \sum_{1\le i <j \le N}[\ee(\jj_{ij}^2;
|\jj_{ij}| > \epsilon\sqrt{N}) + \ee(\jj^{\prime 2}_{ij};
|\jj^\prime_{ij}|>\epsilon \sqrt{N})]  + \epsilon.
\end{eqnarray*}
This shows that the limit of the free energy is the same as that in
the i.i.d.\ standard Gaussian case whenever $\jj_{ij}$'s are
independent with zero mean and unit variance, and satisfy
\begin{eqnarray}\label{mincond1}
& & \forall \epsilon > 0, \ \lim_{N\ra \infty} N^{-2} \sum_{1\le i < j\le
  N} \ee (\jj_{ij}^2; |\jj_{ij}| > \epsilon \sqrt{N}) = 0.
\end{eqnarray}
Note that this is almost exactly condition (\ref{mincond}), the only
difference being that here we do not have terms corresponding to
$i=j$. In particular, it is satisfied when $\jj_{ij}$'s are i.i.d.\
with zero mean and unit variance.

Under the assumption of uniformly bounded third absolute moments,
Corollary \ref{cor1} can be applied to get an explicit error bound of
order $N^{-1/2}$, which is the same as that obtained by Carmona and Hu 
\cite{carmonahu04}.

\section{Ground state of the S-K model}\label{ground}
The ground state in a spin glass model is the configuration which
minimizes the Hamiltonian. With $\beta = 1$ and $h = 0$ for
simplicity, the energy of the ground state is given by
\[
S_N(\jj) = \max_{\bs \in \Sigma^N} \sum_{1\le i < j \le N} \jj_{ij}
\sigma_i \sigma_j.
\]
Guerra and Toninelli \cite{guerra03, guerratoninelli02} proved that 
$N^{-3/2} S_N(\jj)$ converges almost surely and in average to a
deterministic limit if $\jj$ is a collection of standard Gaussian
random variables. It was extended to the case of i.i.d.\ entries
with zero mean, unit variance and finite third moment by Carmona and Hu
\cite{carmonahu04}. We shall show that convergence in probability and 
in average to the same limit would hold if $\jj_{ij}$'s were
independent and satisfied the same condition as in the previous section.

Let $\mf$, $\jj$ and $\jj^\prime$ be as in the previous section, with
$\beta = 1$ and $h = 0$. If we
let $U_N = \max_{\bs} f_{\bs} (\jj)$ and $V_N = \max_{\bs} f_{\bs}
(\jj^\prime)$, then by Theorem 
\ref{maxthm}, for any thrice differentiable $g:\rr \ra \rr$ and any
fixed $K>0$ and $\alpha \ge 1$,  $|\ee g(U_N) - \ee g(V_N)|$ is bounded by
a constant multiple (depending only on $g$) of 
\begin{eqnarray*}
\alpha^{-1} N + \alpha N^{-3} \sum_{i<j}[\ee(\jj_{ij}^2;
|\jj_{ij}| > K) + \ee(\jj^{\prime 2}_{ij}; |\jj^\prime_{ij}|>K)]  +
\alpha^2 N^{-5/2} K.
\end{eqnarray*}
Now choose any $A \ge 1$ and $\epsilon > 0$, and put $\alpha = AN$ and 
$K = \epsilon \sqrt{N}$. Substituting these values in the above
expression, we get
\[
A^{-1} + A  N^{-2}\sum_{i<j}[\ee(\jj_{ij}^2; |\jj_{ij}| > \epsilon
\sqrt{N}) + \ee(\jj^{\prime 2}_{ij};
|\jj^\prime_{ij}|>\epsilon\sqrt{N})] + A^2 \epsilon.
\]
Thus, under condition (\ref{mincond1}) of the previous section,
$\limsup_{N \ra \infty} |\ee g(U_N) - \ee g(V_N)| \le A^{-1} + 
A^2 \epsilon$. This proves the claim, since $A$ and $\epsilon$ are
arbitrary. 

Again, Corollary \ref{cor2} can be applied to obtain an error bound of 
order $N^{-1/6}$ under the assumption of uniformly bounded third
absolute moments.

\section{Proofs}\label{proofs}
{\sc Proof of Theorem \ref{boundthm} } As mentioned before, the
proof is just an easy extension of  
Lindeberg's argument for the classical central limit theorem. Fix
$f$ and $g$ as in the statement of the Theorem. Let $h = 
g \circ f$. Then observe that 
\begin{small}
\begin{eqnarray*}
\partial_i^2 h(\bx) &=& \gp(f(\bx)) \partial_i^2 f(\bx) + \gpp(f(\bx)) 
(\partial_i f(\bx))^2, \\
\partial_i^3 h(\bx) &=& \gp(f(\bx)) \partial_i^3 f(\bx) + 3 \gpp(f(\bx))
\partial_i f (\bx)\partial_i^2 f(\bx) + \gppp(f(\bx)) (\partial_i f(\bx))^3.
\end{eqnarray*}
\end{small}
It follows that for any $i$ and $\bx$, $|\partial_i^2 h(\bx)| \le
C_1 \lambda_2(f)$ and $|\partial_i^3 h(\bx)| \le 6C_2
\lambda_3(f)$, where $C_1 = \|\gp\|_\infty + \|\gpp\|_\infty$ and
$C_2 = \frac{1}{6}\|\gp\|_\infty + \frac{1}{2}\|\gpp\|_\infty +
\frac{1}{6}\|\gppp\|_\infty$.    

Next, for $0 \le i \le n$,  
define $\mathbf{Z}_i := (X_1,\ldots,X_{i-1}, X_i,Y_{i+1},\ldots, Y_n)$ 
and $\mathbf{W}_i := (X_1,\ldots,X_{i-1}, 0, Y_{i+1},\ldots,Y_n)$,
with obvious meanings for $i=0$ and $n$. For $1\le i\le n$, define
\begin{eqnarray*}
R_i &:=& h(\mathbf{Z}_i) - X_i \partial_i h(\mathbf{W}_i) - \frac{1}{2}X_i^2
\partial_i^2 h(\mathbf{W}_i), \\
T_i &:=&  h(\mathbf{Z}_{i-1}) - Y_i \partial_i h(\mathbf{W}_i) -
\frac{1}{2}Y_i^2\partial_i^2 h(\mathbf{W}_i). 
\end{eqnarray*}
By third order Taylor expansion
and the bounds on the third partials of $h$ obtained above, we immediately
see that $|R_i|\le C_2\lambda_3(f)|X_i|^3$ and 
$|T_i|\le C_2\lambda_3(f) |Y_i|^3$. Second order bounds, on the other
hand, imply that $|R_i| \le C_1 \lambda_2(f) |X_i|^2$ and $|T_i| \le C_1
\lambda_2(f)|Y_i|^2$.
Now for each $i$, $X_i$, $Y_i$ and $\mathbf{W}_i$ are independent. Hence 
\begin{eqnarray*}
\ee(X_i \partial_i f(\mathbf{W}_i)) - \ee(Y_i \partial_i
f(\mathbf{W}_i)) &=& \ee(X_i - Y_i) \ee(\partial_i f(\mathbf{W}_i))
\ = \ 0.
\end{eqnarray*}
Similarly, $\ee(X_i^2 \partial_i^2 f(\mathbf{W}_i)) - \ee(Y_i^2
\partial_i^2 f(\mathbf{W}_i)) = 0$.
Combining all these observations we have, for any $K > 0$, 
\begin{eqnarray*}
|\ee g(U) - \ee g(V)| &=& |\sum_{i=1}^n \ee(h(\mathbf{Z}_i) -
h(\mathbf{Z}_{i-1}))| \\
&=& |\sum_{i=1}^n \ee(X_i \partial_i h(\mathbf{W}_i) +
\frac{1}{2}X_i^2 \partial_i^2 h(\mathbf{W}_i) + R_i) \\
& & \ - \ \sum_{i=1}^n \ee(Y_i \partial_i
h(\mathbf{W}_i) + \frac{1}{2}Y_i^2 \partial_i^2 h(\mathbf{W}_i) + T_i)
| \\
&\le& C_1 \lambda_2(f)\sum_{i=1}^n  [\ee (X_i^2; |X_i| > K) + \ee 
(Y_i^2; |Y_i| > K)] \\ 
& & \!\!+ \ C_2 \lambda_3(f) \sum_{i=1}^n  [\ee (|X_i|^3; |X_i| \le K)
+ \ee (|Y_i|^3; |Y_i| \le K)].  
\end{eqnarray*}
The corollary follows by taking $K \ra \infty$.  \hfill $\Box$ \\

{\sc Proof of Theorem \ref{likelihood} } We begin by defining a bunch
of functions. The domains 
will be clear from the definitions. Let
\begin{eqnarray*}
\psi(\bx, f) &:=& e^{\alpha f(\bx)}, \\ 
Z(\bx) &:=& \sum_{f\in \mf} \psi(\bx, f), \\
p(\bx,f) &:=& Z(\bx)^{-1} \psi(\bx, f), \\
a_i(\bx, f) &:=& \alpha \partial_i f(\bx), \\
e_i(\bx) &:=& \sum_{f\in \mf} a_i (\bx, f) p(\bx, f).
\end{eqnarray*}
Note that for any $\bx$, $p(\bx,\cdot)$ is a probability on
$\mf$. This will be widely used without mention in
obtaining the bounds below. Also, note that $F(\bx) =  \alpha^{-1}
\log Z(\bx)$.  

We shall now find bounds on the partial derivatives of several orders
for these functions. Function arguments will be suppressed for
clarity. First, note that clearly from the given expressions,  
\begin{eqnarray}
\partial_i \psi &=& a_i \psi, \label{df} \\
\partial_i Z &=& \sum_{f\in \mf} \partial_i \psi \ = \ \sum_{f\in \mf}
a_i \psi \ = \ Z e_i. \label{dz} 
\end{eqnarray}
Using (\ref{df}) and (\ref{dz}) and the expression for $p$ we get
\begin{equation}\label{dp}
\partial_i p = \frac{Z a_i \psi - Z e_i \psi}{Z^2} = (a_i - e_i) p. 
\end{equation}
Now, directly from the expression for $e_i$ we get
\begin{eqnarray}
\partial_i e_i &=& \sum_{f\in \mf}
(p\partial_i a_i  + a_i \partial_i p), \label{d1e}\\
\partial^2_i e_i &=& \sum_{f\in \mf}
(p \partial^2_i a_i + 2(\partial_i a_i) (\partial_i p) + a_i
\partial^2_i p). \label{d2e}
\end{eqnarray}
Using (\ref{dp}) and (\ref{d1e}) we get
\begin{equation}\label{d2p}
\partial^2_i p = (\partial_i a_i - \partial_i e_i) p + 
(a_i - e_i)^2 p.
\end{equation}
Now for $1\le r\le 3$, let $C_r = \sup\{|\partial_i^r f(\bx)|: 1\le
i\le n, f \in \mf, \ \bx \in I^n\}$. Then note that for any $i$ we
have the uniform bounds    
\begin{equation}\label{abd}
|a_i| \le \alpha C_1, \ |\partial_i a_i| \le 
\alpha C_2, \ |\partial^2_i a_i| \le \alpha C_3
\end{equation}
In the following, we shall freely use the assumption that $\alpha \ge
1$. The first inequality above immediately gives  
\begin{equation}\label{ebd}
|e_i| \le \alpha C_1.
\end{equation}
From (\ref{dp}), (\ref{abd}) and (\ref{ebd}), we get
\begin{equation}\label{dpbd}
|\partial_i p| \le 2\alpha C_1 p.
\end{equation}
Using (\ref{d1e}), (\ref{abd}) and (\ref{dpbd})  we get
\begin{equation}\label{debd}
|\partial_i e_i| \le \alpha^2(C_2 + 2 C_1^2).
\end{equation}
Using (\ref{d2p}), (\ref{abd}), (\ref{ebd}) and (\ref{debd}) we get
\begin{equation}\label{d2pbd}
|\partial^2_i p| \le \alpha^2(2C_2 + 6C_1^2)p.
\end{equation}
Using (\ref{d2e}), (\ref{abd}), (\ref{dpbd}) and
(\ref{d2pbd}) we have
\begin{equation}\label{d2ebd}
|\partial^2_i e_i| \le \alpha^3(C_3 + 6C_1C_2 + 6C_1^3).
\end{equation}
The proof is completed by observing that $\partial_i F = \alpha^{-1}
\partial_i \log Z =  \alpha^{-1} e_i$ 
and using the bounds (\ref{ebd}), (\ref{debd}) and
(\ref{d2ebd}) in Definition \ref{lamdef}.  \hfill $\Box$\\

{\sc Proof of Theorem \ref{maxthm} } For each $\alpha \ge 1$, let
$F_\alpha(\bx) = \alpha^{-1} \log [\sum_{f\in \mf}e^{\alpha
  f(\bx)}]$. Also, let $F(\bx) = \max_{f\in\mf} f(\bx)$. Then we have 
\begin{eqnarray*}
F(\bx) &=& \alpha^{-1} \log [e^{\alpha \max_{f\in \mf} f(\bx)}] \\
&\le& \alpha^{-1} \log [\sum_{f\in \mf} e^{\alpha f(\bx)}] \\
&\le& \alpha^{-1} \log [|\mf| e^{\alpha \max_{f\in \mf} f(\bx)}] 
\end{eqnarray*}
which gives the uniform bound
\[
|F(\bx) - F_\alpha(\bx)| \le \alpha^{-1} \log |\mf|.
\]
Thus, by Theorem \ref{likelihood}, for any $K > 0$,
\begin{eqnarray*}
|\ee g(F(\bbx)) - \ee g(F(\bby))| &\le& 2\|\gp\|_\infty\alpha^{-1}
\log|\mf| + 3\alpha C_1(g) \lambda_2(\mf) T_1(K) \\
& & \ + \ 13\alpha^2 C_2(g)
\lambda_3(\mf) T_2(K)
\end{eqnarray*}
where $T_1(K) = \sum_{i=1}^n [\ee (X_i^2; |X_i| > K) + \ee (Y_i^2;
|Y_i| > K)]$ and $T_2(K) = \sum_{i=1}^n [\ee (|X_i|^3; |X_i| \le K) +
\ee (|Y_i|^3; |Y_i| \le K)]$. If $\gamma <\infty$, then we can let $K
\ra \infty$ and get
\[
|\ee g(F(\bbx)) - \ee g(F(\bby))| \le 2\|\gp\|_\infty\alpha^{-1}
\log|\mf| + 26 \alpha^2 C_2(g) \lambda_3(\mf) \gamma n.
\]
Now choose $\alpha = [(\gamma n\lambda_3(\mf))^{-2/3}(\log|\mf|)^{2/3} +
1]^{1/2}$. Note that $\alpha \ge 1$ and $\alpha^{-1} \le (\gamma
n\lambda_3(\mf))^{1/3} (\log|\mf|)^{-1/3}$. The Corollary follows from 
this. \hfill $\Box$\\
\\
{\bf Acknowledgement. } The author thanks Persi Diaconis for helpful
comments and encouragement, and Erwin Bolthausen for communicating the 
work of Carmona and Hu.

\begin{small}

\end{small}
\end{document}